# A NEW METHOD OF WEIGHT MULTIPLICITIES EVALUATION FOR SEMI-SIMPLE LIE ALGEBRAS

## ANATOLI LOUTSIOUK

Abstract. In most applications of semi-simple Lie groups and algebras representation theory, calculating weight multiplicities is one of the most often used and effort consuming operations. The existing tools were created many years ago by Kostant and Freuedenthal. The celebrated Kostant weight multiplicity formula uses summation over the Weyl group of values of Kostant partition function, and the Freudenthal formula is recurrent. In this paper, a new way for calculating weight multiplicities is presented. The method does not employ the Weyl group, and it is direct, not recurrent. The algorithm realized in accordance with this method is much faster than those realized with the previously employed techniques. Many examples of programs realized in Fortran language are given. They are ready for compilation and execution on desktop and laptop personal computers.

1. **INTRODUCTION**

Given a complex semi-simple Lie algebra G of rank n, let H be a Cartan subalgebra of G. Consider an irreducible finite-dimensional representation T of G in a complex linear space $E^T$ of dimension $m^T$. A non-zero vector v of the space $E^T$ is said to be a weight vector of the representation T if there exists a linear form α on the Cartan subalgebra H so that for all h∈ $H$, T(h)v=α(h)v. The linear form α is called a weight of the representation T. For a given weight α, the subset $E^α$ of the space $E^T$ consisting of all vectors of weight α is a linear subspace called a weight subspace of the space $E^T$, and its dimension is called the multiplicity of the weight α. The space $E^T$ is a direct sum of the weight subspaces, so that $m^T$ is equal to the sum of all the weight multiplicities of the representation T. For a given irreducible finite-dimensional representation T, there exists a unique nonnegative integral weight $p_T$ of multiplicity 1, which is called the highest weight of the representation and which determines the representation up to equivalence. On the other hand, for any non-negative integral linear form β on H there exists a unique up to equivalence irreducible finite-dimensional representation of G for which β is the highest weight. Any integral weight β on H can be identified with a sequence of n integers β = ($β_1$, $β_2$,…,$β_n$), where the integers $β_i$ are called the weight components. Another important mathematical object related to a complex semi-simple Lie algebra G is a finite group $W_G$, which is called the Weyl group of the Lie algebra. For any finite-dimensional representation T of G, the group $W_G$ acts on the set of all weights of representation T. For any weight β of a representation T and for any element w of the group $W_G$ the weight wβ of representation T has the same multiplicity as the weight β. There is one particular finite-dimensional representation of the Lie algebra G, which is called the adjoint representation. The space of this representation is the Lie algebra G itself. The non-zero weights of the adjoint representation are called the roots of the Lie algebra G. The roots are weights of multiplicity 1, and the family of all roots is subdivided into two subfamilies - positive roots and negative roots. There is one element in the dual space of Cartan subalgebra H which plays a special role in the representation theory. This element is equal to the half-sum of all positive roots and is denoted by ρ. The multiplicity of a weight β of an irreducible finite-dimensional representation T with the highest weight α can be evaluated by the following Kostant formula [2]:

Mult(β)=$\sum_{s\in W_G}(\det s)P(s(α + ρ) − (β + ρ))$. (1)

In this formula there is an undefined symbol P. P denotes the Kostant partition function. For any integral vector μ of H', P(μ) equals the number of all representations of μ as a linear combination of positive roots with non-negative integer coefficients. In this nice-looking formula, there are two things that are hard to deal with: the first one is the Weyl group that is needed explicitly and which is very large for a Lie algebra of high rank. Evaluation of Kostant function is also very complicated for an algebra of high rank.

Another formula mostly used for evaluation of weight multiplicities is Freudenthal recurrent formula[1]. This formula does not need Weyl group, and does not need partition function, but it has some computational difficulties, because of a huge amount of computations in the case of high ranks and large weights.

$$((\alpha+\rho,\alpha+\rho) - (\beta+\rho,\beta+\rho))\text{mult}(\beta) = 2 \sum_{\gamma>0} \sum_{k=1}^{\infty} mult(\beta + k\gamma)(\beta + k\gamma, \gamma) \quad (2)$$

In this formula, the first sum is over the set of all positive roots, and ( , ) denotes an inner product on the dual space H' of the linear space H.

Semi-simple Lie algebras are direct sums of simple Lie algebras. Simple Lie algebras are well known. There are four infinite series: $A_n$ (n≥1), $B_n$ (n≥2), $C_n$ (n≥3), $D_n$ (n≥4), and five exceptional algebras: $G_2$, $F_4$, $E_6$, $E_7$, $E_8$. In all cases, the value of index n is equal to the rank of the simple Lie algebra.

## 2. A NEW APPROACH TO THE EVALUATION OF WEIGHT MULTIPLICITIES

Given a complex semi-simple Lie algebra G of rank n, and its Cartan subalgebra H with dual space H'. The roots and weights are elements of H'. In the dual space H' we consider the basis consisting of simple roots, and we identify every positive root and weight with the sequence of its coordinates with respect to this basis. Order the set $\Delta^+$ of all positive roots so that the first n elements in this ordering are the simple roots. Denote by t the number of elements in the set $\Delta^+$. We now introduce one more real linear space $R^t$. This space has dimension t and in it we identify the elements of the canonical basis with the elements of the set $\Delta^+$ of all positive roots, so that the first n vectors are identified with the simple roots. We consider the linear space $R^n$ as a linear subspace of the space $R^t$ generated by the first n vectors of the basis. For an irreducible finite dimensional representation $T=T_p$, with the highest weight $p_T=(p_1,p_2,...,p_k)$, consider the set $W^T$ of all weights of T as a subset of the space $R^n$ as follows: shift the lowest weight vector into the origin of $R^n$, so that the set of weights moves to the subset of the space $\mathbb{R}^n$ with all non-negative components. To each weight realized in this way, we attach its multiplicity.

**Statement 1.** For any irreducible representation T(p) of G there exists a convex polytope F(p) in the space $R^t$ which completely determines the weight multiplicities for T(p) as follows: for any weight μ of the representation

$$\text{mult}(\mu) = LP_p(\mu) \quad (3)$$

where $LP_p$ denotes the local partition function, that is the restriction of Kostant partition function P to the polytope F(p). By restriction of the Kostant partition function to F(p), we mean that only those linear combinations of positive roots with nonnegative integer coefficients which belong to F(p) are counted. The positive roots that are not simple are replaced by their expressions as sums of simple roots.

**Statement 2.** For any complex semi-simple Lie algebra G and for any its irreducible representation T(p) for which all weight components are equal to a nonnegative integer s, the polytope F(p) is the t-dimensional cube of all the vectors $v=(v_1, v_2,...,v_t)$ with $0 \le v_i \le s$ for all I.

**Statement 3.** In the general case of an arbitrary semi-simple Lie algebra G of rank n and an arbitrary irreducible representation T(p) with the highest weight $p=(p_1, p_2,...,p_n)$, the convex polytope F(p) is a subset of the parallelepiped of all the vectors $v=(v_1,...,v_n,v_{n+1},...,v_t)$ with $0 \le v_i \le p_i$ for all i from 1 to n, and $0 \le v_i \le \max(p_1,...,p_n)$ for all i from n+1 to t.

**Statement 4.** For any vector $v=(v_1, ...,v_n)$ with all components $v_i$ satisfying the conditions $0 \le v_i \le p_i$ for all values of i from 1 to n, the vector $\tilde{v}=(v_1,..,v_n,0,...,0)$ belongs to the polytope F(p).

**Statement 5.** If $m=\max(p_1,p_2,...,p_n)$, then for any i>n, the vector $v^i=(0,...,0,...,m,0,..0)$, with i*th* component equal to m and all other components equal to 0, belongs to F(p).

**Statement 6.** Jf a vector v=(v$_1$, v$_2$,....,v$_t$) belongs to F(p), then any its subordinate vector with nonnegative integer components u=(u$_1$,u$_2$,...,u$_t$) with u$_i$ ≤ v$_i$ for all i, also belongs to F(p).

**Remark.** Unfortunately, the amazing simplicity of the polytope F(p) described in Statement 2 for all p with all components of p equal to the same integer does not hold in the case of the highest weight p not belonging to this class; nevertheless for any semi-simple Lie algebra G the principal weights p can be subdivided into subclasses for which the respective polytopes are described by the same formulas with values of p as parameters. For example, in the case of G belonging to the class A$_2$, the number of such subclasses equals 2, in the case of B$_2$ the number of subclasses is equal to 3, for G$_2$ the number of subclasses equals 5.

**Statement 7.** The subclasses of highest weights of which it is said in the Remark above are defined as follows: 1) Consider a set A of n linearly independent positive roots (n is the rank of G) with the property that in the decomposition of any positive root α not belonging to this set with respect to A at least one coordinate with respect to this set is strictly negative. 2) The set of all linear combinations of elements of A with nonnegative coefficients is called a chamber. (Analog of Weyl Chamber). All irreducible representations T(p) of G with highest weights p belonging to the same chamber have polytopes F(p) described by similar formulas.

### 3. CODES OF PROGRAMS IN FORTRAN EVALUATING WEIGHT MULTIPLICITIES FOR IRREDUCIBLE REPRESENTATIONS OF SOME SIMPLE COMPLEX LIE ALGEBRAS

This section contains programs in Fortran language of algorithms realizing calculations on the basis of Statements 1 – 7.

The first program evaluates weight diagrams for all irreducible representations of a Lie algebra of class A$_2$. The highest weight of a representation is given as a vector of 2 nonnegative integers (p,q). The values of p and q must be entered into the program in the third line of the code. Then the program is ready for compilation and execution. The output is done to the screen and to a text file with name 'MULTIPLICITIES.TXT'. The positive roots are stored in the integer array named a. The polytope is stored in the logical array named k. The program has 3 subroutines, which evaluate the polytopes k depending on the type of the highest weight: (p=q, p<q, p>q). In fact, there are only 2 chambers, (p≤q) and (p≥q), the case (p=q) is the common border between the 2 chambers, and so the first subroutine can be excluded from the code.

```
PROGRAM a2multPtoQ        !! Finds weight multiplicities for A2 with (p,q) all cases
IMPLICIT NONE
INTEGER, PARAMETER::p=5,q=3
INTEGER, PARAMETER::n=MAX(p,q)
LOGICAL,DIMENSION(0:p,0:q,0:n)::k=.FALSE.
INTEGER,DIMENSION(0:2*n,0:2*n)::m=0 !multiplicites of weights
INTEGER,DIMENSION(3,2)::a=0 !Positive roots
INTEGER::i1,i2,i3,j1,j2,dim
OPEN(UNIT=1,FILE='MULTIPLICITIES.TXT',STATUS='REPLACE',ACTION='WRITE')
a(1,1)=1; a(2,2)=1; a(3,1)=1; a(3,2)=1 !Nonnegative components of positive roots
IF (p==q) THEN   !Calling subroutines depending on the type of the highest weight.
 CALL a2first(p,q,n,k)
```

```fortran
    ELSE IF (p<q) THEN
     CALL a2second(p,q,n,k)
    ELSE  IF(q<p) THEN
     CALL a3third(p,q,n,k)
    END IF

    DO j1=0,2*n; DO j2=0,2*n
    m(j1,j2)=0

    DO i1=0,p          !Calculating multiplicities
    DO i2=0,q
    DO i3=0,n
    IF(k(i1,i2,i3) .AND. &
            ((i1*a(1,1)+i2*a(2,1)+i3*a(3,1)==j1).AND.(i1*a(1,2)+i2*a(2,2)+i3*a(3,2)==j2)))THEN
    m(j1,j2)=m(j1,j2)+1
    END IF
    END DO; END DO; END DO; END DO; END DO

    DO j1=0,2*n                 !writes multiplicities to screen and to file
    DO j2=0,2*n
     WRITE(*,*)'m(',j1,',',j2,')=',m(j1,j2)
     WRITE(1,*)'m(',j1,',',j2,')=',m(j1,j2)
    END DO; END DO
    dim=0
    DO j1=0,2*n; DO j2=0,2*n
    dim=dim+m(j1,j2)
    END DO; END DO
    WRITE(*,*)'Dimension by summation of weight multiplicities is equal to ',dim
    WRITE(*,*)'Dimension by formula      is equal to ', ((p+1)*(q+1)*(p+q+2))/2
    END PROGRAM a2multPtoQ

    SUBROUTINE a2first(p,q,n,k)            !!p==q
    IMPLICIT NONE
    INTEGER::p,q,n,i1,i2,i3
```

```
LOGICAL,DIMENSION(0:p,0:q,0:n)::k
k=.TRUE.
END SUBROUTINE a2first

SUBROUTINE a2second(p,q,n,k)                    !!p<q
IMPLICIT NONE
INTEGER::p,q,n,i1,i2,i3
LOGICAL,DIMENSION(0:p,0:q,0:n)::k
DO i1=0,p; DO i2=0,p; DO i3=0,n
 k(i1,i2,i3)=.TRUE.
END DO; END DO; END DO
DO i1=0,p; DO i2=p+1,q; DO i3=0,n
k(i1,i2,i3)=((i1+i2+i3)<=p+q)
END DO; END DO; END DO
END SUBROUTINE a2second

SUBROUTINE a2third(p,q,n,k)                     !!q<p
IMPLICIT NONE
INTEGER::p,q,n,i1,i2,i3
LOGICAL,DIMENSION(0:p,0:q,0:n)::k
DO i1=0,q; DO i2=0,q; DO i3=0,n
k(i1,i2,i3)=.TRUE.
END DO; END DO; END DO
DO i1=q+1,p; DO i2=0,q; DO i3=0,n
k(i1,i2,i3)= ((i1+i2+i3)<=p+q)
END DO; END DO; END DO
END SUBROUTINE a2third
```

The next program evaluates weight multiplicities for all weights of an irreducible representation T of the Lie algebra of class $G_2$ with the highest weight (n,n) . For n=9, the dimension of the representation T with the highest weight (9, 9) is equal to 1 million. The program evaluates the weight diagram (multiplicities of all weights of T) on an ordinary notebook personal computer within 2 minutes.

**PROGRAM** g2multforequalweights*! Finds all multiplicities for representations of G2 with highest weight (n,n).*

IMPLICIT NONE

```fortran
INTEGER, PARAMETER::n=9
INTEGER, DIMENSION(0:10*n,0:10*n)::m=0 ! weight multiplicities values
INTEGER,DIMENSION(6,2)::a=0 !Positive roots
INTEGER::i1,i2,i3,i4,i5,i6,j1,j2,dim
OPEN(UNIT=1,FILE='MULTIPLICITIES.TXT',STATUS='REPLACE',ACTION='WRITE')
! Nonnegative components of positive roots.
a(1,1)=1; a(2,2)=1; a(3,1)=1; a(3,2)=1; a(4,1)=2; a(4,2)=1; a(5,1)=3; a(5,2)=1; a(6,1)=3; a(6,2)=2
DO j1=0,10*n; DO j2=0,10*n
m(j1,j2)=0
DO i1=0,n; DO i2=0,n; DO i3=0,n; DO i4=0,n; DO i5=0,n; DO i6=0,n
IF(((i1*a(1,1)+i2*a(2,1)+i3*a(3,1)+i4*a(4,1)+i5*a(5,1)+i6*a(6,1)==j1).AND.(i1*a(1,2)+i2*a(2,2)+i3*a(3,2)+i4*a(4,2)+i5*a(5,2)+i6*a(6,2)==j2))) THEN
m(j1,j2)=m(j1,j2)+1
END IF
END DO; END DO; END DO; END DO; END DO; END DO; END DO; END DO
DO j1=0,10*n; DO j2=0,10*n
WRITE(*,*) m(j1,j2)
WRITE(1,*)'m(',j1,',',j2,')=',m(j1,j2)
END DO; END DO
dim=0
DO j1=0,10*n; DO j2=0,10*n
dim=dim+m(j1,j2)
END DO; END DO
WRITE(*,*) 'Dimension by summation of multiplicities is equal to ', dim
WRITE(*,*) 'Dimension by formula is equal to', ((n+1)*(n+1)*(n+n+2)*(n+2*n+3)*(n+3*n+4)*(2*n+3*n+5))/120
END PROGRAM g2multforequalweights
```

The next program evaluates weight multiplicities for all irreducible representations of a Lie algebra of class $A_4$, with the highest weight of the kind (n,n,n,n) (all weight components of the highest weight equal)

```fortran
PROGRAM a4multequal         !computes to file weight multiplicities for representations of A4 with highest weight (n,n,n,n)
IMPLICIT NONE
INTEGER,PARAMETER::n=3
INTEGER, DIMENSION(0:4*n,0:6*n,0:6*n,0:4*n)::m=0
```

```fortran
INTEGER,DIMENSION(10,4)::a=0 !Positive roots
INTEGER::i1,i2,i3,i4,i5,i6,i7,i8,i9,i10,j1,j2,j3,j4,dim
OPEN(UNIT=1,FILE='MULTIPLICITIES.TXT',STATUS='REPLACE',ACTION='WRITE')
! Nonnegative components of positive roots.
a(1,1)=1; a(2,2)=1; a(3,3)=1; a(4,4)=1; a(5,1)=1; a(5,2)=1; a(6,2)=1; a(6,3)=1; a(7,3)=1; a(7,4)=1
a(8,1)=1; a(8,2)=1; a(8,3)=1; a(9,2)=1;  a(9,3)=1;a(9,4)=1; a(10,1)=1; a(10,2)=1; a(10,3)=1; a(10,4)=1
DO j1=0,4*n; DO j2=0,6*n; DO j3=0,6*n; DO j4=0,4*n
m(j1,j2,j3,j4)=0
DO i1=0,n; DO i2=0,n; DO i3=0,n; DO i4=0,n; DO i5=0,n; DO i6=0,n; DO i7=0,n; DO i8=0,n; DO i9=0,n; DO i10=0,n
IF(((i1*a(1,1)+i2*a(2,1)+i3*a(3,1)+i4*a(4,1)+i5*a(5,1)+i6*a(6,1)+i7*a(7,1)+i8*a(8,1)+i9*a(9,1)+i10*a(10,1)==j1)&
.AND.(i1*a(1,2)+i2*a(2,2)+i3*a(3,2)+i4*a(4,2)+i5*a(5,2)+i6*a(6,2)+i7*a(7,2)+i8*a(8,2)+i9*a(9,2)+i10*a(10,2)==j2)&
.AND.(i1*a(1,3)+i2*a(2,3)+i3*a(3,3)+i4*a(4,3)+i5*a(5,3)+i6*a(6,3)+i7*a(7,3)+i8*a(8,3)+i9*a(9,3)+i10*a(10,3)==j3)&
.AND.(i1*a(1,4)+i2*a(2,4)+i3*a(3,4)+i4*a(4,4)+i5*a(5,4)+i6*a(6,4)+i7*a(7,4)+i8*a(8,4)+i9*a(9,4)+i10*a(10,4)==j4))) THEN
m(j1,j2,j3,j4)=m(j1,j2,j3,j4)+1
END IF
END DO; END DO; END DO; END DO; END DO; END DO; END DO; END DO; END DO; END DO; END DO; END DO; END DO; END DO
DO j1=0,4*n; DO j2=0,6*n; DO j3=0,6*n; DO j4=0,4*n
IF (m(j1,j2,j3,j4)/=0) THEN
WRITE(*,*) m(j1,j2,j3,j4)
WRITE(1,*)'m(',j1,',',j2,',',j3,',',j4,')=',m(j1,j2,j3,j4)
END IF
END DO; END DO; END DO; END DO
dim=0
DO j1=0,4*n; DO j2=0,6*n; DO j3=0,6*n; DO j4=0,4*n
dim=dim+m(j1,j2,j3,j4)
END DO; END DO; END DO; END DO
WRITE(*,*)'Dimension by summation of multiplicities is equal to ',dim
WRITE(*,*)'Dimension by formula is equal to',(n+1)**10
END PROGRAM a4multequal
```

The next program evaluates weight multiplicity for a particular weight (not the weight diagram, as in other programs presented here) of an irreducible representation T(p) of the Lie algebra of class $F_4$, with the highest weight p=(n, n, n, n)

**PROGRAM** f4multequalweight!!    *Finds multiplicities of weights for repesentations of F4 with highest weight (n,n,n,n)*

IMPLICIT NONE

INTEGER,PARAMETER::n=2

INTEGER::m=0 ! *multiplicity*

INTEGER,DIMENSION(24,4)::a=0 ! *Positive roots.*

INTEGER::i1,i2,i3,i4,i5,i6,i7,i8,i9,i10,i11,i12,i13,i14,i15,i16,i17,i18,i19,i20,i21,i22,i23,i24

INTEGER:: j1, j2, j3, j4, dim

WRITE(*,*)'Enter weight, 4 nonnegative integers' !*The components of the weight whose multiplicity you want to evaluate*

READ(*,*) j1, j2, j3, j4

! *Nonnegative components of positive roots.*

a(1,1)=1; a(2,2)=1; a(3,3)=1; a(4,4)=1; a(5,1)=1; a(5,2)=1; a(6,2)=1; a(6,3)=1; a(7,3)=1; a(7,4)=1

a(8,1)=1; a(8,2)=1; a(8,3)=1; a(9,2)=1; a(9,3)=1; a(9,4)=1; a(10,1)=1; a(10,2)=1; a(10,3)=1; a(10,4)=1

a(11,2)=2; a(11,3)=1; a(12,2)=2; a(12,3)=1; a(12,4)=1; a(13,1)=1; a(13,2)=2; a(13,3)=1; a(14,2)=2; a(14,3)=2

a(14,4)=1; a(15,1)=1; a(15,2)=2; a(15,3)=1; a(15,4)=1; a(16,1)=2; a(16,2)=2; a(16,3)=1; a(17,1)=1; a(17,2)=2

a(17,3)=2; a(17,4)=1; a(18,1)=2; a(18,2)=2; a(18,3)=1; a(18,4)=1; a(19,1)=1; a(19,2)=3; a(19,3)=2; a(19,4)=1

a(20,1)=2; a(20,2)=2; a(20,3)=2; a(20,4)=1; a(21,1)=2; a(21,2)=3; a(21,3)=2; a(21,4)=1; a(22,1)=2; a(22,2)=4

a(22,3)=2; a(22,4)=1; a(23,1)=2; a(23,2)=4; a(23,3)=3; a(23,4)=1; a(24,1)=2; a(24,2)=4; a(24,3)=3; a(24,4)=2

DO i1=0,n;DO i2=0,n;DO i3=0,n;DO i4=0,n;DO i5=0,n;DO i6=0,n;DO i7=0,n;DO i8=0,n;DO i9=0,n;DO i10=0,n;DO i11=0,n

DO i12=0,n;DO i13=0,n;DO i14=0,n;DO i15=0,n;DO i16=0,n;DO i17=0,n;DO i18=0,n;DO i19=0,n;DO i20=0,n;DO i21=0,n

DO i22=0,n;DO i23=0,n;DO i24=0,n

IF(((i1*a(1,1)+i2*a(2,1)+i3*a(3,1)+i4*a(4,1)+i5*a(5,1)+i6*a(6,1)&

+i7*a(7,1)+i8*a(8,1)+i9*a(9,1)+i10*a(10,1)+i11*a(11,1)+i12*a(12,1)&

+i13*a(13,1)+i14*a(14,1)+i15*a(15,1)+i16*a(16,1)+i17*a(17,1)+i18*a(18,1)&

+i19*a(19,1)+i20*a(20,1)+i21*a(21,1)+i22*a(22,1)+i23*a(23,1)+i24*a(24,1)==j1)&

.AND.(i1*a(1,2)+i2*a(2,2)+i3*a(3,2)+i4*a(4,2)+i5*a(5,2)+i6*a(6,2)&

```
        +i7*a(7,2)+i8*a(8,2)+i9*a(9,2)+i10*a(10,2)+i11*a(11,2)+i12*a(12,2)&
        +i13*a(13,2)+i14*a(14,2)+i15*a(15,2)+i16*a(16,2)+i17*a(17,2)+i18*a(18,2)&
        +i19*a(19,2)+i20*a(20,2)+i21*a(21,2)+i22*a(22,2)+i23*a(23,2)+i24*a(24,2)==j2)&
        .AND.(i1*a(1,3)+i2*a(2,3)+i3*a(3,3)+i4*a(4,3)+i5*a(5,3)+i6*a(6,3)&
        +i7*a(7,3)+i8*a(8,3)+i9*a(9,3)+i10*a(10,3)+i11*a(11,3)+i12*a(12,3)&
        +i13*a(13,3)+i14*a(14,3)+i15*a(15,3)+i16*a(16,3)+i17*a(17,3)+i18*a(18,3)&
        +i19*a(19,3)+i20*a(20,3)+i21*a(21,3)+i22*a(22,3)+i23*a(23,3)+i24*a(24,3)==j3)&
        .AND.((i1*a(1,4)+i2*a(2,4)+i3*a(3,4)+i4*a(4,4)+i5*a(5,4)+i6*a(6,4)&
        +i7*a(7,4)+i8*a(8,4)+i9*a(9,4)+i10*a(10,4)+i11*a(11,4)+i12*a(12,4)&
        +i13*a(13,4)+i14*a(14,4)+i15*a(15,4)+i16*a(16,4)+i17*a(17,4)+i18*a(18,4)&
        +i19*a(19,4)+i20*a(20,4)+i21*a(21,4)+i22*a(22,4)+i23*a(23,4)+i24*a(24,4)==j4)))  THEN
m=m+1
END IF
END DO; END DO; END DO; END DO;END DO; END DO; END DO; END DO;END DO; END DO; END DO; END DO;END DO; END DO
END DO; END DO;END DO; END DO; END DO; END DO;END DO; END DO; END DO; END DO
WRITE(*,*)'The multiplicity is =',m
END PROGRAM f4multequalweight
```

The next program finds all multiplicities for any irreducible representation of a Lie algebra of class $B_2$. The program has 4 subroutines, though the number of chambers is equal to 3. The first subroutine can be excluded, because the case (p=q) is a border case between 2 chambers, and can be treated by other subroutines.

```
PROGRAM b2multPtoQ          !!(p,q) all cases for B2
IMPLICIT NONE
INTEGER,PARAMETER::p=10,q=12
INTEGER,PARAMETER::n=MAX(p,q)
LOGICAL,DIMENSION(0:p,0:q,0:n,0:n)::k=.FALSE.   !The polytope.
INTEGER,DIMENSION(0:4*n,0:4*n)::m=0
INTEGER,DIMENSION(4,2)::a=0 !Positive roots.
INTEGER::i1,i2,i3,i4,j1,j2,dim
OPEN(UNIT=1,FILE='MULTIPLICITIES.TXT',STATUS='REPLACE',ACTION='WRITE')
! Nonnegative components of positive roots.
a(1,1)=1;      a(2,2)=1; a(3,1)=1; a(3,2)=1; a(4,1)=1; a(4,2)=2
IF (p==q) THEN
```

```
    CALL b2first(p,q,n,k)
ELSE IF (p<q) THEN
    CALL b2second(p,q,n,k)
ELSE IF((q<p).AND.((2*q)>=p))THEN
    CALL b2third(p,q,n,k)
ELSE IF((q<p).AND.((2*q)<p)) THEN
    CALL b2fourth(p,q,n,k)
END IF

DO j1=0,4*n; DO j2=0,4*n
m(j1,j2)=0
        !Calculates multiplicities
DO i1=0,p; DO i2=0,q; DO i3=0,n; DO i4=0,n
IF(k(i1,i2,i3,i4) .AND.  ((i1*a(1,1)+i2*a(2,1)+i3*a(3,1)+i4*a(4,1)==j1).AND.
        (i1*a(1,2)+i2*a(2,2)+i3*a(3,2)+i4*a(4,2)==j2)))              THEN
m(j1,j2)=m(j1,j2)+1
END IF; END DO; END DO; END DO; END DO; END DO; END DO

DO j1=0,4*n              !writes multiplicities to screen and to file
DO j2=0,4*n
 WRITE(*,*)'m(',j1,',',j2,')=',m(j1,j2)
 WRITE(1,*)'m(',j1,',',j2,')=',m(j1,j2)
END DO; END DO
dim=0
DO j1=0,4*n; DO j2=0,4*n
dim=dim+m(j1,j2)
END DO; END DO
WRITE(*,*)'Dimension by summation is equal to ',dim
WRITE(*,*)'Dimension by formula    is equal to ', ((p+1)*(q+1)*(p+q+2)*(2*p+q+3))/6
END PROGRAM b2multPtoQ

SUBROUTINE b2first(p,q,n,k)                 !!p==q
IMPLICIT NONE
INTEGER::p,q,n,i1,i2,i3,i4
```

LOGICAL,DIMENSION(0:p,0:q,0:n,0:n)::k

k=.TRUE.

**END SUBROUTINE** b2first

**SUBROUTINE** b2second(p,q,n,k)  !!*p<q*

IMPLICIT NONE

INTEGER::p,q,n,i1,i2,i3,i4

LOGICAL,DIMENSION(0:p,0:q,0:n,0:n)::k

DO i1=0,p; DO i2=0,p; DO i3=0,n; DO i4=0,n

k(i1,i2,i3,i4)= (((i2+i3+i4)<=p+q).OR.((i3)<=p))

END DO; END DO; END DO; END DO

DO i1=0,p; DO i2=p+1,q; DO i3=0,n; DO i4=0,n

k(i1,i2,i3,i4)=(((i2+i3+i4)<=(p+q)).AND. ((i1+i2+i3)<=(p+q)))

END DO; END DO; END DO; END DO

**END SUBROUTINE** b2second

**SUBROUTINE** b2third(p,q,n,k)  !!*q<p and 2\*q>=p*

IMPLICIT NONE

INTEGER::p,q,n,i1,i2,i3,i4

LOGICAL,DIMENSION(0:p,0:q,0:n,0:n)::k

DO i1=0,q; DO i2=0,2*q-p; DO i3=0,n; DO i4=0,n

k(i1,i2,i3,i4)=.TRUE.

END DO; END DO; END DO; END DO

DO i1=0,q; DO i2=2*q-p+1,q; DO i3=0,n; DO i4=0,n

IF (((i2+i3+i4)<=p+q).OR.((i1+i2+i3)<=2*q)) THEN

k(i1,i2,i3,i4)=.TRUE.

END IF; END DO; END DO; END DO; END DO

DO i1=q+1,p; DO i2=0,q; DO i3=0,n; DO i4=0,n

IF (((i1+i2+i3)<=p+q)) THEN

k(i1,i2,i3,i4)=.TRUE.

END IF; END DO; END DO; END DO; END DO

**END SUBROUTINE** b2third

**SUBROUTINE** b2fourth(p,q,n,k)            !!$(q<p)$ and$(2*q<p)$

IMPLICIT NONE

INTEGER::p,q,n,i1,i2,i3,i4

LOGICAL,DIMENSION(0:p,0:q,0:n,0:n)::k

DO i1=0,q; DO i2=0,q; DO i3=0,n; DO i4=0,n

k(i1,i2,i3,i4)=(((i2+i3+i4)<=p+q).OR.((i3)<=q))

END DO; END DO; END DO; END DO

DO i1=q+1,p-q-1; DO i2=0,q; DO i3=0,n; DO i4=0,n

k(i1,i2,i3,i4)=(((i2+i3+i4)<=(p+q)).AND. ((i1+i2+i3)<=(p+q)).OR.(i3<=q))

END DO; END DO; END DO; END DO

DO i1=p-q,p; DO i2=0,q; DO i3=0,n; DO i4=0,n

k(i1,i2,i3,i4)=(i3<=q)

END DO; END DO; END DO; END DO

**END SUBROUTINE** b2fourth

The following program evaluates all multiplicities of all weights of irreducible representations of a Lie algebra of class $G_2$ with the highest weight (0, q), with q any positive integer.

**PROGRAM** g2multwith0toq !!Calculates weight multiplicities with weight (0,q) for G2

IMPLICIT NONE

INTEGER,PARAMETER::p=0,q=9

!!        *Finds all multiplicities for representations of G2 with highest weight (0,q)*

INTEGER,PARAMETER::n=MAX(p,q)

INTEGER, DIMENSION(0:10*n,0:10*n)::m=0

INTEGER,DIMENSION(6,2)::a=0 !*Positive roots*

INTEGER::i1,i2,i3,i4,i5,i6,j1,j2,dim

LOGICAL,DIMENSION(0:p,0:q,0:n,0:n,0:n,0:n)::k=.FALSE.

OPEN(UNIT=1,FILE='MULTIPLICITIES.TXT',STATUS='REPLACE',ACTION='WRITE')

! *Nonnegative components of positive roots.*

a(1,1)=1; a(2,2)=1; a(3,1)=1; a(3,2)=1; a(4,1)=1; a(4,2)=2; a(5,1)=1; a(5,2)=3; a(6,1)=2; a(6,2)=3

DO i1=0,p; DO i2=0,q; DO i3=0,n; DO i4=0,n; DO i5=0,n; DO i6=0,n

k(i1,i2,i3,i4,i5,i6)= (((i3+i4+i5+i6)<=q).AND.((i2+i3+i4+i5)<=q))        !!!*polytope*

END DO; END DO; END DO; END DO; END DO; END DO

DO j1=0,10*n; DO j2=0,10*n

m(j1,j2)=0

DO i1=0,p; DO i2=0,q; DO i3=0,n; DO i4=0,n; DO i5=0,n; DO i6=0,n

IF (k(i1,i2,i3,i4,i5,i6).AND.(((i1*a(1,1)+i2*a(2,1)+i3*a(3,1)+i4*a(4,1)+i5*a(5,1)+i6*a(6,1)==j1).AND.(i1*a(1,2)+i2*a(2,2)+i3*a(3,2)+i4*a(4,2)+i5*a(5,2)+i6*a(6,2)==j2)))) THEN

m(j1,j2)=m(j1,j2)+1

END IF

END DO; END DO; END DO; END DO; END DO; END DO; END DO; END DO

DO j1=0,10*n; DO j2=0,10*n

IF (m(j1,j2)>0) THEN

WRITE(*,*)'m(',j1,',',j2,')=',m(j1,j2)

WRITE(1,*)'m(',j1,',',j2,')=',m(j1,j2)

END IF

END DO; END DO

dim=0

DO j1=0,10*n; DO j2=0,10*n

dim=dim+m(j1,j2)

END DO; END DO

WRITE(*,*)'Dimension is equal to ',dim

WRITE(*,*)'Dimension by formula is equal to',((q+1)*(p+1)*(q+p+2)*(q+2*p+3)*(q+3*p+4)*(2*q+3*p+5))/120

**END PROGRAM** g2multwith0toq

The following program evaluates all multiplicities of all weights of irreducible representations of a Lie algebra of class $G_2$ with the highest weight (p, 0), with p any positive integer,

**PROGRAM** g2multpto0         !!! calculates multiplicities for representations of $G_2$ with highest weight (p,0)

IMPLICIT NONE

INTEGER,PARAMETER::p=3,q=0

!!      Finds all multiplicities for representations of G2 with highest weight (p,0)

INTEGER,PARAMETER::n=MAX(p,q)

```fortran
INTEGER, DIMENSION(0:10*n,0:10*n)::m=0
INTEGER,DIMENSION(6,2)::a=0 !Positive roots.
INTEGER::i1,i2,i3,i4,i5,i6,j1,j2,dim
LOGICAL,DIMENSION(0:p,0:q,0:n,0:n,0:n,0:n)::k=.FALSE.
OPEN(UNIT=1,FILE='MULTIPLICITIES.TXT',STATUS='REPLACE',ACTION='WRITE')
OPEN(UNIT=2,FILE='k.txt',status='REPLACE',ACTION='WRITE')
! Nonnegative components of positive roots.
a(1,1)=1; a(2,2)=1; a(3,1)=1; a(3,2)=1; a(4,1)=1; a(4,2)=2; a(5,1)=1; a(5,2)=3; a(6,1)=2; a(6,2)=3
DO i1=0,p; DO i2=0,q; DO i3=0,n; DO i4=0,n; DO i5=0,n; DO i6=0,n
k(i1,i2,i3,i4,i5,i6)=(((i4+i5)<=p).AND.((i1+i3)<=p).AND.((i3+i4)<=p).AND.((i4+i6)<=p))   !!!polytope
IF ((k(i1,i2,i3,i4,i5,i6)).eqv.(.TRUE.))    THEN
WRITE(*,*)'k(',i1,',',i2,',',i3,',',i4,',',i5,',',i6,')=',k(i1,i2,i3,i4,i5,i6)
END IF
END DO; END DO; END DO; END DO; END DO; END DO

DO j1=0,10*n; DO j2=0,10*n
m(j1,j2)=0

DO i1=0,p; DO i2=0,q; DO i3=0,n; DO i4=0,n; DO i5=0,n; DO i6=0,n
IF (k(i1,i2,i3,i4,i5,i6).AND.(((i1*a(1,1)+i2*a(2,1)+i3*a(3,1)+i4*a(4,1)+i5*a(5,1)+i6*a(6,1)==j1).AND.(i1*a(1,2)+i2*a(2,2)+i3*a(3,2)+i4*a(4,2)+i5*a(5,2)+i6*a(6,2)==j2)))) THEN
m(j1,j2)=m(j1,j2)+1
END IF
END DO; END DO; END DO; END DO; END DO; END DO; END DO; END DO
DO j1=0,10*n; DO j2=0,10*n
IF (m(j1,j2)>0) THEN
WRITE(*,*)'m(',j1,',',j2,')=',m(j1,j2)
WRITE(1,*)'m(',j1,',',j2,')=',m(j1,j2)
END IF
END DO; END DO
dim=0
DO j1=0,10*n; DO j2=0,10*n
dim=dim+m(j1,j2)
```

END DO; END DO

WRITE(*,*)'Dimension by summation ofmultiplicities is equal to ',dim

WRITE(*,*)'Dimension by formula is equal to',((q+1)*(p+1)*(q+p+2)*(q+2*p+3)*(q+3*p+4)*(2*q+3*p+5))/120

**END PROGRAM** g2multpto0

DEPT.CIVIL ENG., KING MONGKUT UNIVERSITY OF TECHNOLOGY THOMBURY, THAILAND

*E-mail address:* loutsiouk.ana@kmutt.ac.th